\documentclass[onefignum,onetabnum]{siamart171218}

\usepackage{pgf}

\usepackage{multirow}
\usepackage{amsfonts}
\usepackage{graphicx}
\usepackage{epstopdf}
\usepackage{algorithmic}
\usepackage{subcaption}
\usepackage{tikz}
\usepackage{enumitem}
\usetikzlibrary{shapes,arrows}
\usetikzlibrary{cd}
\usetikzlibrary{decorations}
\usepackage{tikz-cd}
\usetikzlibrary{decorations.pathmorphing}
\usetikzlibrary{positioning}
\usetikzlibrary{shapes.multipart}
\ifpdf
  \DeclareGraphicsExtensions{.eps,.pdf,.png,.jpg}
\else
  \DeclareGraphicsExtensions{.eps}
\fi

\newsiamremark{remark}{Remark}
\newsiamremark{hypothesis}{Hypothesis}
\newsiamremark{assumption}{Assumption}
\crefname{hypothesis}{Hypothesis}{Hypotheses}
\newsiamthm{claim}{Claim}

\headers{Polynomial approximation for moment matching}{Carlos Doebeli et al.}

\title{A polynomial approximation scheme for nonlinear model reduction by moment matching\thanks{Submitted to the editors December 16, 2024; accepted for publication (in revised form) September 15, 2026.}}

\author{Carlos Doebeli\thanks{Department of Applied Mathematics, Imperial College London (\email{c.doebeli22@imperial.ac.uk}, \email{d.kalise-balza@imperial.ac.uk}).}
\and Alessandro Astolfi\thanks{Computer, Electrical and Mathematical Science and Engineering Division, King Abdullah University of Science and Technology (KAUST), Thuwal, 23955, Saudi Arabia (\email{alessandro.astolfi@kaust.edu.sa}).}
\and Dante Kalise\footnotemark[2]
\and Alessio Moreschini\footnotemark[4]
\and Giordano Scarciotti\thanks{Department of Electrical and Electronic Engineering, Imperial College London (\email{a.moreschini@imperial.ac.uk}, \email{g.scarciotti@imperial.ac.uk}).}
\and Joel Simard \thanks{Formerly at the Department of Electrical and Electronic Engineering, Imperial College London. Presently independent researcher (\email{joeldsimard@gmail.com}).}
}

\usepackage{amsopn}

\makeatletter
\newcommand*{\addFileDependency}[1]{
  \typeout{(#1)}
  \@addtofilelist{#1}
  \IfFileExists{#1}{}{\typeout{No file #1.}}
}
\makeatother

\ifpdf
\hypersetup{
  pdftitle={A polynomial approximation scheme for nonlinear model reduction by moment matching},
  pdfauthor={C. Doebeli, A. Astolfi, D. Kalise, A. Moreschini, G. Scarciotti, J. Simard}
}
\fi

\begin{document}

\maketitle
\begin{abstract}
We propose a procedure for the numerical approximation of invariance equations arising in the moment matching technique associated with reduced-order modeling of high-dimensional dynamical systems. The Galerkin residual method is employed to find an approximate solution to the invariance equation using a Newton iteration on the coefficients of a monomial basis expansion of the solution. These solutions to the invariance equations can then be used to construct reduced-order models. We assess the ability of the method to solve the invariance PDE system as well as to achieve moment matching and recover the steady-state behaviour of nonlinear systems with state dimension of order 1000 driven by linear and nonlinear signal generators.
\end{abstract}



\section{Introduction}

Complex large-scale dynamical systems are ubiquitous in control engineering, providing a natural framework for modeling multi-agent systems, energy networks, and coupled electrical systems \cite{khalil}. The high dimensionality of these models, often on the order of thousands or millions of states, introduces significant computational challenges, a phenomenon known as the ``curse of dimensionality''~\cite{bellman}. Reduced-order modeling aims to alleviate this complexity while preserving important system characteristics. It finds application across many fields, including mechanical and electrical systems and weather forecasting \cite{MOR_app, large_scale_dynam}. Broadly speaking, model order reduction methods are categorized as being based upon the singular value decomposition, such as Balanced Truncation~\cite{kramer2022nonlinear, kramer2023nonlinear, scherpen1993balancing} and Proper Orthogonal Decomposition~\cite{Kunisch1999ControlOT, Kunisch2008, peraire}, or upon Krylov projection/interpolation, such as moment matching~\cite{astolfi_mm, astolfi_mm_nonlinear}; see~\cite{antoulas_review} for a comprehensive overview. Other approaches include energy-based methods \cite{gray_hilbert, gray_energy}, reduction around a limit cycle \cite{gray_manifolds}, and transfer function interpolation for polynomial systems \cite{benner}.

In this paper, we consider high-dimensional systems with a state
$x(t) \in \mathbb{R}^n$, input $u(t) \in \mathbb{R}^m$, and output $y(t) \in \mathbb{R}^p$,
governed by
\begin{equation*}
    \dot{x} = f(x,u), \qquad y = h(x).
\end{equation*}
In particular, we are interested in the case where the input $u$ is generated
by an exogenous system with state $\omega(t) \in \mathbb{R}^d$, described by
\begin{equation*}
    \dot{\omega} = s(\omega), \qquad v = \ell(\omega).
\end{equation*}
Under certain assumptions, the interconnected system obtained by setting $u = v$
admits a mapping $x = \pi(\omega)$ satisfying a system of first-order nonlinear
partial differential equations (PDEs) known as the invariance equation, which characterizes an invariant manifold
embedding the interconnected system's dynamics. The goal of this work is to
achieve nonlinear model reduction by moment matching by solving this invariance
equation, yielding a reduced system with state $r(t) \in \mathbb{R}^d$ governed by
\begin{equation*}
    \dot{r} = \bar{f}(r,u), \qquad y_r = \bar{h}(r),
\end{equation*}
whose output $y_r$ approximates the steady-state behaviour of the full-order
system under the same input. The precise mathematical setting is given in
\Cref{sec:mm}.

For recent surveys on moment matching see \textit{e.g.} \cite{rafiq2022, scarciotti2024interconnection}. For linear systems, moment matching is in one-to-one relation with the solution to a Sylvester equation induced by the interconnection of the system and the signal generator \cite{gallivan}; with additional stability assumptions the invariant manifold is attractive and links to the steady-state response \cite{isidori2008steady}. This relationship was extended to nonlinear systems in \cite{isidori1995nonlinear}, associating invariance with a ``Sylvester-like'' system of nonlinear coupled PDEs whose (locally) attractive solutions generalize the notion of steady-state response. Building on this, \cite{astolfi_mm} generalized moment matching to nonlinear systems via the invariant manifold induced by the interconnection. Further developments include two-sided nonlinear moment matching \cite{ionescu2016nonlinear}, data-driven estimation via basis function fitting \cite{scarciotti2017data, moreschini2024a}, parametrization of all differential-algebraic interpolants \cite{simard2023parameterization, simard2023parameterizationb}, and a closed-loop framework relaxing stability requirements \cite{moreschini2024closed}. The notion of invariance also arises beyond model reduction, in output regulation \cite{byrnes2003limit} and in immersion and invariance control design \cite{astolfi_ii}.

The main limitation of the interconnection-based approach is the need to numerically solve a large-scale system of coupled nonlinear PDEs: for an $n$-dimensional system and a $d$-dimensional signal generator, this amounts to solving $n$ PDEs over $\mathbb{R}^d$. In this paper, we propose a numerical scheme based on Newton's method in conjunction with a global polynomial ansatz for each PDE solution. This approach, first explored for Hamilton-Jacobi-Bellman PDEs in \cite{beard_galerkin, beard_thesis} and later extended to high-dimensional settings in \cite{kalise, kalise2019robust, AKK21}, mitigates the curse of dimensionality at moderate computational cost. We focus on low-dimensional signal generators with $d = 2$ and high-dimensional dynamics with $n$ up to 1000, a scale for which traditional grid-based PDE methods are impractical. The use of orthogonal polynomials for weighted residual methods in moment matching was explored in \cite{faedo_scarciotti_astolfi} for linear signal generators and low-dimensional dynamics.

\subsection{Contribution}

This paper focuses on the development of computational methods for nonlinear moment matching. The main contributions are as follows:

\begin{itemize}
    \item We develop and implement a Galerkin method to find approximate solutions to the invariance PDE system and use these to construct reduced-order models achieving moment matching.

    \item The framework is tailored to nonlinear systems driven by low-dimensional signal generators, exploiting the interconnection structure to construct very low-dimensional reduced-order models that accurately reproduce steady-state behaviour and outperform standard model-order reduction (MOR) techniques at the same reduced order.

    \item This paper represents the first computational method explicitly developed for approximating the invariance equation in nonlinear moment matching. The method's ability to handle systems with thousands of states is a key novel contribution.

    \item The method's efficacy is demonstrated on systems of dimension up to 1000, systems driven by non-polynomial signal generator nonlinearities, and interconnected systems with implicitly defined nonlinear signal generators.

    \item The method applies beyond model order reduction, to output regulation and immersion and invariance problems.
\end{itemize}

\subsection{Organization} The paper is organized as follows. In \Cref{sec:mm} we give an overview of the problem setup and we show how moment matching is defined and how the invariance PDE system arises. \Cref{sec:galerkin} introduces a Galerkin type method for the numerical computation of an approximate solution $\pi(\cdot)$ to the invariance PDE system. \Cref{sec:results} reports different numerical experiments assessing the accuracy of the proposed method and its effectiveness in generating a reduced-order model. \Cref{sec:discussion} discusses the results further and puts them in a wider context.

\section{Moment Matching for Nonlinear Dynamical Systems} \label{sec:mm}

We consider dynamical systems with a state $x(t) \in \mathbb{R}^n$, an
input $u(t) \in \mathbb{R}^m$, and an output $y(t) \in \mathbb{R}^p$, that can
be represented in the form
\begin{equation}\label{eq:system}
\begin{aligned}
    \dot{x} & = f(x,u),  \\ 
    y & = h(x). \\
\end{aligned}
\end{equation}
We assume that the mappings $f : \mathcal{U} \times \mathcal{V} \to \mathbb{R}^n$
and $h : \mathcal{U} \to \mathbb{R}^p$ are $C^2$, defined on open
neighbourhoods $\mathcal{U} \subset \mathbb{R}^n$ and $\mathcal{V} \subset \mathbb{R}^m$
of the origin.
We assume that system \eqref{eq:system} has an equilibrium and, without loss of generality, that this equilibrium is satisfied at the origin, that is, $f(0,0) = 0$ and $h(0) = 0$.


In order to introduce the notion of moment for a nonlinear system of the form \eqref{eq:system},  we consider the interconnection of \eqref{eq:system} with an exogenous system, referred to as the \textit{signal generator}, with a state $\omega(t) \in \mathbb{R}^d$ which is described by equations of the form
\begin{equation}\label{eq:sig_system}
\begin{aligned}
    \dot{\omega} & = s(\omega), \\
    v & = \ell(\omega), \\
\end{aligned}
\end{equation} 
where $v(t) \in \mathbb{R}^{m}$ is the output of the signal generator system, and the mappings $s : \mathbb{R}^d \to \mathbb{R}^d$ and $\ell : \mathbb{R}^d \to \mathbb{R}^m$ are $C^2$. We also assume that \eqref{eq:sig_system} has an equilibrium at the origin.

We denote by $\Xi_t^s : \mathbb{R} \times \mathbb{R}^d \to \mathbb{R}^d$ the flow generated by the vector field $s$, that is, the smooth function with the property that
\begin{equation}
    \frac{\partial}{\partial t} \Xi_t^s(\omega_0) = s(\Xi_t^s(\omega_0)),  \quad \Xi_0^s(\omega_0) = \omega_0.
\end{equation} 
Here, we introduce the notion of observability and accessibility, using the definitions introduced in \cite{moreschini2024closed}, and based on \cite{hermann} and \cite{sussmann}:

\begin{definition}[Observability]
    The system~\eqref{eq:system} is said to be locally observable at the point $x_0$ if for every open neighbourhood $U$ of $x_0$, initial conditions $x_a(0) = x_0 \neq x_b(0) \in U$, and input $u(\cdot) \in L^{\infty}\left([0,\infty); \mathbb{R}^m\right)$ for which the corresponding trajectories $x_a(\cdot), x_b(\cdot)$ remain in $U$, the output trajectories satisfy $h(x_a(t)) \neq h(x_b(t))$ for some $t > 0$, \textit{i.e.} the output trajectories are not identical.

    If we relax the condition and say that this holds when $U$ is the entire space $\mathcal{U}$, then the system~\eqref{eq:system} is said to be observable at $x_0$. If a system is (locally) observable at every $x_0$, then it is said to be (locally) observable.    
\end{definition}

\begin{definition}[Accessibility]
Let $\mathcal{U} \subset \mathbb{R}^n$ be a set which contains some open ball with center $x = 0$. Let $\mathcal{R}(x(0),T)$ be the set (reachable set) containing all points $\overline{x}$ for which there exists an input $u(\cdot) \in L^{\infty}\left([0,\infty); \mathbb{R}^m\right)$ such that the evolution of \eqref{eq:system} from $x(0)\in\mathcal{U}$ satisfies $x(t)\in\mathcal{U}$ for $0 \leq t \leq T$ with $x(T) = \overline{x}$. The system~\eqref{eq:system} is said to be (locally) accessible if for all $x(0)\in\mathcal{U}$ the set $\bigcup_{t\leq T}\mathcal{R}(x(0),t)$ contains a non-empty open subset of $\mathcal{U}$ for all $T>0$.
\end{definition}

Throughout this article, we will assume that the system \eqref{eq:system} is locally observable and locally accessible. We will also assume that the signal generator \eqref{eq:sig_system} is locally observable. From here, it is useful to establish definitions of Poisson stability and neutral stability of a system in connection to the signal generator.

\begin{definition} [Poisson stability]
    A point $\omega_0 \in \mathbb{R}^d$ is said to be Poisson stable for system \eqref{eq:sig_system} if the flow $\Xi_t^s(\omega_0)$ generated by the vector field $s(\omega)$ is defined for all $t \in \mathbb{R}$ and, for each neighbourhood $U^0$ about $\omega_0$ and each $T > 0$, there exists $t_1, t_2 \in \mathbb{R}$ with $t_1 > T$ and $t_2 < -T$, such that the flow satisfies $\Xi_{t_1}^s(\omega_0) \in U^0$ and $\Xi_{t_2}^s(\omega_0) \in U^0$.
\end{definition}

\begin{definition} [Neutral stability]
    A dynamical system is said to be neutrally stable if $\omega = 0$ is a stable equilibrium, and there is an open neighbourhood $W^0$ around the origin in which every point is Poisson stable. 
\end{definition}

Clearly, if the signal generator is neutrally stable on a domain $W^0$, then $\omega$ is persistent in time and does not decay to zero as time tends towards infinity, for any initial condition $\omega_0 \in W^0$. Furthermore, this implies that the linear approximation of the system around the origin $S = \left. \left[\frac{\partial s} {\partial \omega}\right] \right|_{\substack{\omega=0}}$ has purely imaginary eigenvalues, and that these eigenvalues are semi-simple \cite[Chapter 8]{isidori1995nonlinear}. 

\subsection{Notion of Moment}

This paper focuses on the interconnected system \eqref{eq:system}-\eqref{eq:sig_system} with the interconnection equation $u = v$. The interconnected system has a state-space representation given by

\begin{equation} \label{eq:interconnected}
\begin{aligned}
    \dot{\omega} & = s(\omega), \\ 
    \dot{x} & = f(x, \ell(\omega)), \\
    y & = h(x), \\
    \omega(0) & = \omega_0, \\
    x(0) & = x_0.
\end{aligned}
\end{equation}

The notion of the time-domain moment for a nonlinear dynamical system comes from analyzing the steady-state output of the interconnected system given by \eqref{eq:interconnected} and makes use of the centre manifold theory. In order to help in our analysis of this interconnected dynamical system, we first define the notion of an invariant manifold and a centre manifold, which are introduced in more detail in \cite[Appendix B]{isidori1995nonlinear}, \cite{carr}, \cite[Chapter 8]{khalil}.

\begin{definition} [Invariant Manifold]
    Consider a system
    \begin{equation} \label{eq:center_sys}
        \dot{\xi} = \eta(\xi),
    \end{equation}
    for $\xi(t) \in \mathbb{R}^n$. A set $\mathcal{M}$ is said to be an invariant manifold for \eqref{eq:center_sys} if $\xi(0) \in \mathcal{M} \Rightarrow \xi(t) \in \mathcal{M}$ for all $t \in \mathbb{R}$.
\end{definition}

\begin{definition} [Centre Manifold]
    Now consider a system of the form
    \begin{equation} \label{eq:center_connected}
        \begin{aligned}
            \dot{\omega} & = S \omega + \tilde{s}(\omega, x) \\
            \dot{x} & = A x + \tilde{f}(\omega, x),
        \end{aligned}
    \end{equation}
    where $\omega(t) \in \mathbb{R}^d, x(t) \in \mathbb{R}^n$, $S$ is a $(d \times d)$ matrix having all eigenvalues with zero real part, $A$ is an $(n \times n)$ matrix having all eigenvalues with negative real part, and the functions $\tilde{s}$ and $\tilde{f}$ are $C^2$ with the functions and their first order derivatives vanishing at the origin. If $\mathcal{M} = \{(\omega,x) : x = \pi(\omega)\}$, for some function $\pi : \mathbb{R}^d \to \mathbb{R}^n$, is an invariant manifold for \eqref{eq:center_connected}, it is also called a centre manifold if it is smooth with $\pi(0) = 0$ and $\frac{\partial \pi}{\partial \omega}(0) = 0$. 

    Moreover, the centre manifold $\mathcal{M}$ characterized by the mapping $x = \pi(\omega)$ is said to be a locally attractive centre manifold for \eqref{eq:center_connected} if there exists a neighbourhood $U^0 \subset \mathbb{R}^d \times \mathbb{R}^n$ of $(0,0)$ and constants $C \geq 1, \ \alpha > 0$ such that if $(\omega(0), x(0)) \in U^0$, then 
    \begin{equation} \label{eq:attract}
    \|x(t) - \pi(\omega(t)) \| \leq C e^{- \alpha t} \|x(0) - \pi(\omega(0)) \|,
    \end{equation}    

    for all $t \geq 0$.

\end{definition}

Under certain assumptions for the signal generator \eqref{eq:sig_system} and the system \eqref{eq:system}, we can assume that a locally attractive centre manifold exists between the two systems. We therefore make the following assumptions about the functions:
\begin{enumerate}[label=({A\arabic*}), topsep=5pt, itemsep=5pt]
    \item The signal generator \eqref{eq:sig_system} is neutrally stable. \label{ass:stab}
    
    \item The system \eqref{eq:system} is locally exponentially stable in the first approximation about the equilibrium at the origin. That is, the matrix $A = \left. \left[\frac{\partial f} {\partial x}\right] \right|_{\substack{x = 0 \\ u = 0}}$ has only eigenvalues with negative real part. \label{ass:f}
\end{enumerate}
The following result is a consequence of \cite[Proposition 8.1.1]{isidori1995nonlinear}:

\begin{proposition}\label{prop}
Consider the interconnected system \eqref{eq:interconnected}. If Assumptions \ref{ass:stab} and \ref{ass:f} hold, then we are guaranteed the existence of a mapping $\pi : W^0 \subset \mathbb{R}^d \to \mathbb{R}^n$, defined on a neighbourhood $W^0$ about the origin, that characterizes a locally attractive centre manifold for \eqref{eq:interconnected}. This mapping $\pi$ satisfies, for all $\omega \in W^0$,
\begin{equation}\label{eq:pde}
\begin{aligned}
    \frac{\partial \pi}{\partial \omega} (\omega) s(\omega) & = f(\pi(\omega), \ell(\omega)), \\
    \pi(0) & = 0 .
\end{aligned}
\end{equation}

For each flow $\Xi_t^s(\omega_0)$ with initial condition $\omega_0 \in W^0$, the mapping $\pi$ solving \eqref{eq:pde} produces a well-defined steady-state response in the interconnected system~\eqref{eq:interconnected}, given by
\begin{equation}
    y_{\mathrm{ss}}(t) = h(\pi(\Xi_t^s(\omega_0))).
\end{equation}
\end{proposition} 
\begin{remark}
Note that in the case where the system \eqref{eq:system} and the signal generator \eqref{eq:sig_system} are both linear, with $s(\omega) = S \omega, \ell(\omega) = L \omega$ and $f(x,u) = Ax + Bu$, the invariance equation \eqref{eq:pde} simplifies to the well-known Sylvester equation
\begin{equation}
    \Pi S = A \Pi + B L
\end{equation}
where the matrix $\Pi$ represents the invariant mapping $x = \Pi \omega$.
\end{remark}

Because this centre manifold is locally attractive, \eqref{eq:attract} implies that for all initial conditions $\omega_0 \in W^0$, if the initial condition $x_0$ of the system is sufficiently close to the centre manifold, the output $y(t;\omega_0,x_0)$ generated by initial conditions $(\omega_0, x_0)$ satisfies
\begin{equation}
    \lim_{t \to \infty} \| y(t;\omega_0, x_0) - h(\pi(\Xi_t^s(\omega_0)))\| = 0.
\end{equation}

The steady-state response of the system is therefore defined by the invariant mapping between $\omega$ and $x$, which allows us to define what we mean by a system's moment.

\begin{definition} [Moment]
    Assume that \eqref{eq:pde} has a solution $\pi$ for the interconnected system \eqref{eq:interconnected}. Then the time-domain moment of the interconnected system at $(s, \ell)$ is given by the function $h(\pi(\omega(\cdot)))$. 
\end{definition}

Assumptions \ref{ass:stab} and \ref{ass:f} guarantee that \eqref{eq:pde} has a solution, and that the notion of moment is in one-to-one relation with the steady-state output of the interconnected system \eqref{eq:interconnected} for initial conditions starting near the centre manifold $\mathcal{M}$. In the rest of the paper, we shall assume that assumptions \ref{ass:stab} and \ref{ass:f} hold.

\subsection{Moment-Based Reduced-Order Modeling} \label{sec:rom}

With the help of our previously defined notion of moment, we would like to use moment matching to develop a reduced-order model that is able to predict the long-term behaviour of the interconnected system, given a known signal generator input. The reduced-order model takes the form 
\begin{equation} \label{eq:reduced}
\begin{aligned}
    \dot{r} & = \bar{f}(r, u) \\
    y_r & = \bar{h}(r)
\end{aligned}
\end{equation}
with $r(t) \in \mathbb{R}^v$, with $v \ll n$, $u(t) \in \mathbb{R}^m$, and where $\bar{f}$ is a $C^2$ mapping $\bar{f} : \mathbb{R}^v \times \mathbb{R}^m \to \mathbb{R}^v$ satisfying $\bar{f}(0,0) = 0$. The output $y_r(t) \in \mathbb{R}^p$ is measured by the $C^2$ mapping $\bar{h} : \mathbb{R}^v \to \mathbb{R}^p$.

The goal is to construct a system of dimension $v \ll n$ whose output $y_r(\cdot)$ matches the steady-state behaviour of $y(\cdot)$ under the same signal generator input. We will use the following definition of a reduced-order model based on \cite{astolfi_mm_nonlinear}:

\begin{definition}
    The system described by \eqref{eq:reduced} is a model at $s(\omega)$ of the system \eqref{eq:system} if it has the same moment at $s(\omega)$ as \eqref{eq:system}. In this case, it is said to match the moment of the system at $s(\omega)$. Furthermore, system \eqref{eq:reduced} is a reduced-order model of system \eqref{eq:system} if it is a model at $s(\omega)$ of \eqref{eq:system} and $v < n$. 
\end{definition}

If the reduced-order model satisfies assumption \ref{ass:f}, then we know that there exists a mapping $p(\cdot)$ satisfying the equation
\begin{equation} \label{eq:ROM_PDE}
\begin{aligned}
    \frac{\partial p}{\partial \omega}(\omega) s(\omega) & = \bar{f}(p(\omega), \ell(\omega)) \\
    p(0) & = 0.
\end{aligned}
\end{equation}

Under these assumptions, the steady-state output of \eqref{eq:reduced} from an initial condition $\omega_0 \in W^0$, with $r_0$ sufficiently close to the centre manifold, satisfies
\begin{equation}\label{eq:reduced_ss}
    \lim_{t \to \infty} | y_r(t; \omega_0, r_0) - \bar{h}(p(\Xi_t^s(\omega_0))) | = 0.
\end{equation}

Therefore, the reduced-order model \eqref{eq:reduced} \textit{matches the moments} of \eqref{eq:system} at $s(\omega)$ if the mapping that solves \eqref{eq:ROM_PDE} also satisfies
\begin{equation} \label{eq:matching}
    h(\pi(\omega)) = \bar{h}(p(\omega)).
\end{equation}
Figure \ref{fig:tikz} depicts how moment matching corresponds to the matching of the steady-state responses of the full-order model and the reduced-order model. 

\tikzstyle{block} = [draw, fill=white, rectangle, 
    minimum height=3em, minimum width=6em]
\tikzstyle{sum} = [draw, fill=white, circle, node distance=1cm]
\tikzstyle{input} = [coordinate]
\tikzstyle{output} = [coordinate]
\tikzstyle{pinstyle} = [pin edge={to-,thin,black}]

\begin{figure}[t!]
\begin{center}
\centering
\begin{tikzpicture}
    \node [block] (sg1) {$\dot{\omega} = s(\omega)$};

    \node [block, right of=sg1, node distance=4.75cm] (state1) {\begin{tabular}{c} $\dot{x} = f(x,u)$ \\ $y = h(x) $\end{tabular}};

    \draw [->] (sg1) -- node[above] {$u = \ell(\omega)$} (state1);
    \node [block, right of=state1, node distance=4.75cm] (output1) {$y_{\mathrm{ss}}(t) = h(\pi(\omega(t)))$};

    \draw [->] (state1) -- node [above] {$y(t)$}(output1);

    \node [text width = 4cm, below right=0.5cm and -3.5cm of output1] (moments) {$\Big\Updownarrow $ (Moment Matching) };

    \node [block, below=2cm of output1] (output2) {$y_{\mathrm{ss}}(t) = \bar{h}(p(\omega(t)))$};

    \node [block, left of=output2, node distance = 4.75cm] (state2) {\begin{tabular}{c} $\dot{r} = \bar{f}(r,u)$ \\ $y_r = \bar{h}(r) $\end{tabular}};

    \node [block, left of=state2, node distance=4.75cm] (sg2) {$\dot{\omega} = s(\omega)$};

    \draw [->] (sg2) -- node[above] {$u = \ell(\omega)$} (state2);

    \draw [->] (state2) -- node [above] {$y_r(t)$}(output2);

\end{tikzpicture}

\caption{Diagrammatic illustration of the concepts of moment matching in reduced-order modeling. The signal generator generates a control input to the high-order model and the reduced-order model, and moment-matching corresponds to the exact matching of their steady-state responses.} 
\label{fig:tikz}
\end{center}
\end{figure}

In order to construct a reduced-order model, we need to choose mappings $\bar{f}$ and $\bar{h}$ such that the equilibrium $(0,0)$ of $\bar{f}$ is locally exponentially stable in the first approximation, and such that \eqref{eq:ROM_PDE} and \eqref{eq:matching} hold. We note that choosing $p(\omega) = \omega$ and $\bar{h}(\omega) = h(\pi(\omega))$ allows the reduced-order model to satisfy \eqref{eq:matching}, and as a consequence simplifies \eqref{eq:ROM_PDE} to
\begin{equation}\label{eq:simplified_pde_rom}
    s(\omega) = \bar{f}(\omega, \ell(\omega)).
\end{equation}
In particular, it is shown in \cite{astolfi_mm} that the choice of $p(\omega) = \omega$ gives us a family of reduced-order models
\begin{equation} \label{eq:constructed_rom}
\begin{aligned}
    \dot{r} & = s(r) - \bar{g}(r, \ell(r)) + \bar{g}(r, u) \\ 
    y_r & = h(\pi(r))
\end{aligned}
\end{equation}
where the mapping $\bar{g}$ is a free parameter, which can be chosen to ensure the local exponential stability about the origin of the autonomous system having dynamics given by
\begin{equation} \label{eq:rom_stab}
    r \mapsto s(r) - \bar{g}(r, \ell(r)).
\end{equation}

A suitable choice of $\bar{g}$ can therefore ensure that Assumption \ref{ass:f} holds, and that the reduced-order model given by \eqref{eq:constructed_rom} achieves moment matching. 

It is convenient to have $\bar{g}$ take the form $\bar{g}(r, \ell(r)) = g(r) \ell(r)$ for some $C^2$ mapping $g : \mathbb{R}^d \to \mathbb{R}^{d \times m}$. In the general case, one suitable mapping can be found by computing
\begin{equation}
\begin{aligned}
    S = \left[\frac{\partial s(r)}{\partial r} \right] \bigg | _{r = 0}, && L =  \left[\frac{\partial \ell(r)}{\partial r} \right] \bigg | _{r = 0}
\end{aligned}
\end{equation}
and choosing a mapping $g(r) = G$ for some real-valued matrix $G$ such that $S - GL$ only has eigenvalues with negative real part, which can be accomplished if $(S,L)$ is detectable. The linearization of \eqref{eq:rom_stab} then simply becomes $S - GL$, and the assumption on the eigenvalues of $S - GL$ ensures that the system is locally exponentially stable. 

\section{Approximate solutions via Galerkin expansion} \label{sec:galerkin}

The main difficulty in constructing a reduced-order model by moment matching lies in obtaining a suitable mapping $\pi(\cdot)$ satisfying the PDE system
\begin{equation}\label{eq:pde_sec3}
\begin{aligned}
    \frac{\partial \pi}{\partial \omega} (\omega) s(\omega) & = f(\pi(\omega), \ell(\omega)), \\
    \pi(0) & = 0 .
\end{aligned}
\end{equation}
A numerical scheme for approximating the solution to \eqref{eq:pde_sec3} is developed in this section. 

\subsection{Polynomial expansion of the invariant mapping}

Proposition \ref{prop} only guarantees the existence of a solution to system \eqref{eq:pde_sec3} in a neighbourhood of the origin, hence we restrict our analysis to a compact domain $\Omega \subset \mathbb{R}^{d}$ containing the origin. We denote the $i$-th component of the $n$-dimensional mapping $\pi$ by $\pi_i$, so $\pi(\omega) = \left( \pi_1(\omega), \pi_2(\omega), \dots, \pi_n(\omega) \right)^\top$. For each mapping $\pi_i : \mathbb{R}^d \to \mathbb{R}$, $i = 1, \dots, n$, we consider an expansion $\pi^N_i$ of the form
\begin{equation} \label{eq:pi_gal}
    \pi^N_i(\omega; \mathbf{c}_i) = \sum_{k=1}^{N} c_{i,k} \phi_k(\omega) = \Phi^N (\omega) \mathbf{c}_i 
\end{equation}
where $\Phi^N = (\phi_1, \phi_2, \dots, \phi_N)$, and where $\{\phi_k\}_{k=1}^{N}$ are the first $N$ elements of a family $\{\phi_k\}_{k=1}^{\infty} \subset C^{\infty}(\Omega)$ that forms a complete system of basis functions in $H^1(\Omega)$.
Note that we choose $H^1(\Omega)$ to ensure that  $\frac{\partial \phi_k}{\partial\omega} \in L^2(\Omega)$. Here, $\mathbf{c}_i = (c_{i,1}, \dots, c_{i,k}, \dots, c_{i,N})^\top \allowbreak \in \mathbb{R}^N$ denotes the vector of coefficients associated with the basis functions $\phi_k$ for the $i$-th component of $\pi$. For this problem, we are interested in globally supported ansatz functions, such as monomials or orthogonal polynomials. The Galerkin residual method finds a weak solution to the invariant manifold $\pi$ by projecting the mapping onto the space of functions $\{\phi_k\}$ and ensuring that the approximation error of the equation is orthogonal to every $\phi_j$ in the basis. The method of weighted residuals presented here and the resulting residual equations shown later are similar to those presented in \cite{faedo_scarciotti_astolfi} and \cite{kalise}.

The components of $\mathbf{c}_i$ in \eqref{eq:pi_gal} for $i = 1, \dots, n$ are found by setting the projection of the error in \eqref{eq:pde_sec3} to zero for all $\omega \in \Omega$ and imposing the Galerkin residual equation
\begin{equation} \label{eq:galerkin_base}
\begin{aligned}
    \left \langle \frac{\partial \pi_i}{\partial \omega} s - f_i(\pi, \ell), \phi_j \right \rangle _{\Omega} = 0, && \text{for } j = 1, \dots, N
\end{aligned}
\end{equation}
where the inner product $\langle\cdot, \cdot \rangle_{\Omega}$ is defined for $f,g \in L^2(\Omega)$ as
\begin{equation} \label{eq:inner}
    \langle f, g \rangle _{\Omega} = \int_{\Omega} f(\omega) g(\omega) d\omega .
\end{equation}

We note that $\frac{\partial \pi_i}{\partial \omega} (\cdot) s(\cdot) - f_i(\pi(\cdot), \ell(\cdot)) \in L^2(\Omega)$ because $\pi$ is smooth by the definition of the centre manifold, $s$ and $f$ are assumed to be at least $C^2$, and $\Omega$ is compact.

For the sake of simplicity, based on the notation used in \cite{beard_thesis}, with a slight abuse of notation we define the following extension of the inner product for vector-valued functions:

\begin{definition} \label{def:inner_prod}
    If $\eta : \Omega \subset \mathbb{R}^N \to \mathbb{R}$ is a real-valued function, then we define 
    \begin{equation}
        \langle \eta, \Phi^N \rangle_{\Omega} = \left( \left \langle \eta, \phi_1 \right \rangle_{\Omega}, \dots, \left \langle \eta, \phi_N \right \rangle_{\Omega} \right)^\top .
    \end{equation}
    If $\eta : \Omega \subset \mathbb{R}^N \to \mathbb{R}^N$ is a real-valued function, then we define 
    \begin{equation} \label{eq:mat_notation}
        \langle \eta, \Phi^N \rangle_{\Omega} = \begin{bmatrix}
            \langle \eta_1, \phi_1 \rangle _{\Omega}& \dots &  \langle \eta_N, \phi_1 \rangle _{\Omega} \\ 
            \vdots & & \vdots \\
            \langle \eta_1, \phi_N \rangle _{\Omega}& \dots &  \langle \eta_N, \phi_N \rangle _{\Omega}
        \end{bmatrix} .
    \end{equation}
\end{definition}

In order to solve the nonlinear system \eqref{eq:galerkin_base}, we introduce the high-dimensional function $F(\cdot) : \mathbb{R}^{Nn} \to \mathbb{R}^{Nn}$, 
\begin{equation} \label{eq:bigF}
    F(\mathbf{c}) := \left( F_1(\mathbf{c}), \dots, F_i(\mathbf{c}), \dots, F_n(\mathbf{c}) \right)^\top
\end{equation}
where each $F_i$ is itself a multi-dimensional function $F_i : \mathbb{R}^{Nn} \to \mathbb{R}^N$, given by
\begin{align} \label{eq:F_i_expanded}
    F_i(\mathbf{c}) &= \begin{bmatrix}
        \left\langle \frac{\partial \pi_i^N}{\partial \omega} s - f_i(\pi^N,\ell), \phi_1 \right\rangle_{\Omega}
        \\
        \vdots
        \\
        \left\langle \frac{\partial \pi_i^N}{\partial \omega} s - f_i(\pi^N,\ell), \phi_N \right\rangle_{\Omega}
    \end{bmatrix} .
\end{align}

By expanding $\pi^N$ in this function, we express each $F_i$ as a sum of linear and nonlinear components on $\mathbf{c}_i$. Using our notation for inner products, we can represent each $F_i(\cdot)$ as

\begin{equation} \label{eq:bigF_i}
\begin{aligned}
    F_{i}(\mathbf{c}) =  \left \langle \frac{\partial \Phi^N}{\partial \omega} s, \Phi^N \right \rangle _{\Omega} \mathbf{c}_i -\left \langle f_i\left( \left[ \Phi^N \mathbf{c}_1, \dots, \Phi^N \mathbf{c}_n \right]^\top, \ell\right), \Phi^N \right \rangle _{\Omega} .
\end{aligned}
\end{equation}

The linear term can be written as a matrix-vector product $\mathbf{A} \mathbf{c}_{i}$, where $\mathbf{A}$ is defined as 
\begin{equation} \label{eq:A_matrix}
\begin{aligned}
    \mathbf{A} & = \left \langle \frac{\partial \Phi^N}{\partial \omega} s, \Phi^N \right \rangle_{\Omega} .
\end{aligned}
\end{equation}

In order to deal with the second term in the right hand side of the expression, for convenience, we can define the function $G : \mathbb{R}^{Nn} \to \mathbb{R}^{Nn}$ as $G(\mathbf{c}) = \left( G_1(\mathbf{c}), \dots, G_n(\mathbf{c})\right)^\top$ with
\begin{equation} \label{eq:G_func}
    G_{i}(\mathbf{c}) = -\left \langle f_i\left( \left[ \Phi^N \mathbf{c}_1, \dots, \Phi^N \mathbf{c}_n \right]^\top, \ell\right), \Phi^N \right \rangle _{\Omega} .
\end{equation}

Using this notation, we can then write $F_i$ as 
\begin{equation} \label{eq:overall_F_i}
    F_i(\mathbf{c}) = \mathbf{A} \mathbf{c}_i + G_i(\mathbf{c})
\end{equation}
and the overall function as
\begin{equation} \label{eq:overall_F}
F(\mathbf{c}) = \begin{bmatrix}
        \mathbf{A} \mathbf{c}_1 \\ \vdots \\ \mathbf{A} \mathbf{c}_n
    \end{bmatrix} + G(\mathbf{c}) .
\end{equation}

\begin{remark} \label{remark:f}
    If the function $f$ determining the system's dynamics is a polynomial of sufficiently low degree, for example if it is linear, quadratic, or cubic in $\pi$, the functions $G_i$ can be represented with matrices and tensors similar to $\mathbf{A}$. In problems where this is the case, it will be useful to construct a matrix or tensor representation of the function in order to accelerate both assembly and iterative computations. 
\end{remark}

In order to find the coefficients $\mathbf{c}$ satisfying $F(\mathbf{c}) = 0$, we employ a Newton iteration, which iterates until $\|F(\mathbf{c})\|_1$ is below a given tolerance. To carry out this iteration, we also need to find the Jacobian of the Galerkin residual function. The Jacobian is defined as $JF : \mathbb{R}^{Nn} \to \mathbb{R}^{Nn \times Nn}$, with $JF(\mathbf{c}) = \frac{\partial F}{\partial \mathbf{c}}(\mathbf{c})$. Based on our formulation in \eqref{eq:bigF_i}, we represent it as 
\begin{equation} \label{eq:JF}
    JF(\mathbf{c}) = \begin{bmatrix}
        JF_{1,1}(\mathbf{c}) & \dots & JF_{1,n}(\mathbf{c}) \\
        \vdots & & \vdots \\
        JF_{n,1}(\mathbf{c}) & \dots & JF_{n,n}(\mathbf{c})
    \end{bmatrix}
\end{equation}
where the corresponding entries are given by
\begin{equation} \label{eq:JF_i,j}
    JF_{i,j}(\mathbf{c}) = \delta_{ij} \mathbf{A} + \frac{\partial G_i}{\partial \mathbf{c}_j}(\mathbf{c}).
\end{equation}

Here, $\delta_{ij}$ denotes the Kronecker delta acting on the indices $i$ and $j$. The functions $\frac{\partial G_i}{\partial \mathbf{c}_j}$ can be computed using the chain rule:
\begin{equation} \label{eq:dG}
    \frac{\partial G_i}{\partial \mathbf{c}_j}(\mathbf{c}) = - \left \langle \frac{\partial f_i(\pi^N, \ell)}{\partial \mathbf{c}_j} \Phi^N,  \Phi^N \right \rangle _{\Omega} \bigg | _{\pi^N = \Phi^N \mathbf{c}}.
\end{equation}

Based on the expression of $F(\mathbf{c})$ in \eqref{eq:bigF_i} and of $JF(\mathbf{c})$ in \eqref{eq:JF_i,j}, we calculate the Newton increment $\Delta \mathbf{c}^k$ by solving

\begin{equation}
    JF(\mathbf{c}^k) \Delta \mathbf{c}^k = F(\mathbf{c}^k),
\end{equation}
and the Newton update is then

\begin{equation}
    \mathbf{c}^{k+1} = \mathbf{c}^k - \Delta\mathbf{c}^k.
\end{equation}

\begin{remark}
    Note that in some cases, depending on the system dynamics under study, the Jacobian $JF$ may become singular. In these scenarios, we solve for the Newton increment $\Delta \mathbf{c}^k$ in the least-squares sense.

\end{remark}

\subsection{Polynomial basis}

The construction of the polynomial basis for our approximation largely follows the procedure used in \cite[Section 4.1]{kalise}. In order to construct the basis functions to be used in the Galerkin approximation over the $d$-dimensional space $\Omega$, we first select a maximum degree $M$ and construct a basis in a single scalar variable $\omega$ up to the degree $M$. We denote this basis by $\varphi_M : \mathbb{R} \to \mathbb{R}^M$. In our case, we will use the monomial basis 
\begin{equation} \label{eq:monomial_basis}
    \varphi_M(\omega) = \left(\omega, \omega^2, \dots, \omega^M \right) .
\end{equation}
Note that, since our invariance equation \eqref{eq:pde_sec3} specified that $\pi(0) = 0$, we do not require a constant term in our basis. In this paper we are using a monomial basis but in general, it is possible to use the same generation method based on a different chosen one-dimensional basis, for example orthogonal bases such as Legendre or Chebyshev polynomials up to degree $M$. 

In order to generate the multi-dimensional basis, we use the tensor product of the basis in each coordinate and eliminate those terms whose total degree is higher than the prescribed maximum degree $M$. The basis is then given by
\begin{equation} \label{eq:multi_basis}
    \Phi^N = \left \{ \phi \in \bigotimes_{j=1}^{d} \varphi_M(\omega_j), \text{ where } \deg(\phi) \leq M \right \}.
\end{equation}
In our case, where the one-dimensional basis is given by the monomials up to degree $M$, the multi-dimensional basis admits a convenient explicit description. Let $\Lambda_M = \{ \nu \in \mathbb{N}_0^d : 1 \leq |\nu| \leq M \}$, with $|\nu| = \sum_{j=1}^d \nu_j$, denote the set of multi-indices of total degree at most $M$ (the zero multi-index is excluded as no constant term is required), so that the number of basis functions is $N = |\Lambda_M| = \binom{M+d}{M} - 1$. We fix an enumeration $i \mapsto \nu(i)$ of $\Lambda_M$, e.g.\ in graded lexicographic order, and write $\nu_{i,j} = \nu(i)_j$; the Galerkin equations are independent of this choice. Each basis polynomial is then the separable function
\begin{equation} \label{eq:basis}
    \phi_i(\omega) = \prod_{j=1}^d \omega_j^{\nu_{i,j}},
    \qquad \nu(i) \in \Lambda_M .
\end{equation}

\subsection{Integration}

We are concerned with the assembly of $\mathbf{A}$, $G_i$, and $\frac{\partial G_i}{\partial \mathbf{c}_j}$. Each component of $\mathbf{A} = \left \langle \frac{\partial \Phi^N}{\partial \omega} s, \Phi^N \right \rangle_{\Omega}$ is given by
\begin{equation} \label{eq:A}
    \mathbf{A}_{i,j} = \left \langle \frac{\partial \phi_j}{\partial \omega} s, \phi_i \right \rangle_{\Omega} = \sum_{k=1}^{d} \int_{\Omega} \frac{\partial\phi_j} {\partial\omega_k} (\omega) s_k(\omega) \phi_i(\omega) d\omega,
\end{equation}
where $s(\omega) = \left( s_1(\omega), \dots, s_d(\omega)\right)^\top$. 
Each component of $G_i(\mathbf{c}) = - \left \langle f_i(\pi^N, \ell), \Phi^N\right \rangle _{\Omega}$ is represented as
\begin{equation} \label{eq:G}
    G_i(\mathbf{c})_{j} = - \left \langle f_i(\pi^N, \ell), \phi_j \right \rangle _{\Omega} = - \int_{\Omega} f_i(\pi^N(\omega ; \mathbf{c}), \ell(\omega)) \phi_j(\omega) d\omega,
\end{equation}
and is approximated using a suitable quadrature rule, for example Gauss-Legendre quadrature in $d$ dimensions over $\Omega$. The entries of $\frac{\partial G_i}{\partial \mathbf{c}_j} = - \left \langle \frac{\partial f_i(\pi^N, \ell)}{\partial \mathbf{c}_j} \Phi^N, \Phi^N\right \rangle _ {\Omega}$ are computed analogously,
\begin{equation} \label{eq:dG_2}
\begin{aligned}
    \left[ \frac{\partial G_i}{\partial \mathbf{c}_j} \right]_{k,l} & = - \left \langle \frac{\partial f_i(\pi^N, \ell)}{\partial \mathbf{c}_j} \phi_l, \phi_k \right \rangle_{\Omega} = - \int_{\Omega} \frac{\partial f_i(\pi^N(\omega ; \mathbf{c}), \ell(\omega))}{\partial \mathbf{c}_j} \phi_l(\omega) \phi_k(\omega) d\omega,
\end{aligned}
\end{equation}
and are likewise approximated by a suitable quadrature rule.

\subsection{Computational complexity and implementation}

The bulk of the computation comes from integration over $d$ dimensions to evaluate
the function $F$ and its Jacobian $JF$ as well as the computation of the Newton
update $\Delta\mathbf{c}^k$. 

The cost of tensorized quadrature scales as $O(q^d)$ for $q$ nodes in each of the dimensions $d$ and thus grows exponentially in $d$, unless $f$ and $s$ are separable or have exploitable structure, as discussed in Remark~\ref{remark:f}. In our setting $d$ is the dimension of the signal generator and is moderate in all examples below, so tensorized Gauss--Legendre quadrature remains practical. For substantially larger $d$, the exponential cost can be mitigated by sparse-grid (Smolyak) quadrature \cite{bungartz_griebel} or quasi-Monte Carlo rules \cite{dick_kuo_sloan}. Plain Monte Carlo integration, with its dimension-independent but slow $O(K^{-1/2})$ convergence, becomes competitive with tensorized quadrature only once $d$ exceeds roughly $8$--$10$ for smooth integrands, since the cost of tensor quadrature at fixed accuracy grows like $(\log1/\varepsilon)^d$ while that of Monte Carlo is $O(\varepsilon^{-2})$
independently of $d$.

When the integrals admit an analytical solution, evaluating $F$ has complexity
$O(Nn)$ and evaluating $JF$ has complexity $O(N^2n^2)$. Solving $JF(\mathbf{c}^k)
\Delta \mathbf{c}^k = F(\mathbf{c}^k)$ via a direct or least-squares solver has cubic
complexity in $Nn$, though in many cases $JF$ is sparse, which significantly
reduces this cost. Overall, the computational cost depends on the system dimension
$n$ and combinatorially on the signal generator dimension $d$ and maximum degree $M$.

\section{Numerical results} \label{sec:results}

Different numerical experiments are performed to assess the approximation quality
of the Galerkin method and the performance of the resulting reduced-order models.
The tests are carried out over rectangular domains $\Omega$ of varying sizes, for
different nonlinear dynamics and signal generators, and the algorithm is employed
with a tolerance of $10^{-7}$ in the $L^1$ norm of $F$. In each case, the
	solver is initialized with $\mathbf{c}^0 = \mathbf{0}$, with all coefficients of $\pi^N$
	set to zero.

The quality of the Galerkin approximation is measured via the residual error
function $R = (R_1, \dots, R_n)^\top$, where for each $i = 1, \dots, n$ the
residual at a specified $\omega$ is given by
\begin{align} \label{eq:error}
	R_i(\mathbf{c}, \omega) = \frac{\partial \pi_i^N}{\partial \omega} (\omega ; \mathbf{c}) s(\omega) - f_i(\pi^N(\omega ; \mathbf{c}), \ell(\omega)).
\end{align}
To quantify the approximation quality globally, we compute the $L_2$ norm of each
$R_i$ over a chosen test domain $W$,
\begin{equation}
	\lVert R_i \rVert _{(2, W)} = \sqrt{ \int_{W} |R_i(\mathbf{c}, \omega)|^2 d\omega},
\end{equation}
and aggregate these into a single metric weighted by the relative magnitude of the
coefficients $\mathbf{c}_i$:
\begin{equation}
    \mathcal{E}_R(W) = \frac{\sum_{i=1}^n \lVert \mathbf{c}_i \rVert_2 \lVert R_i \rVert _{(2, W)} } {\sum_{i=1}^n \lVert \mathbf{c}_i \rVert _2}.
\end{equation}
When a reduced-order model can be constructed, its performance is assessed by
	comparing steady-state outputs. Given time series $y_k = y(t_k)$ and $y_{r,k} =
	y_r(t_k)$ sampled at times $\{t_k\}_{k=1}^{K}$, we identify a time $T_{\mathrm{ss}}$
	beyond which both models have reached steady state, and set $k_s = \min \{k : t_k
	\geq T_{\mathrm{ss}} \}$. Denoting by $\frac{1}{2} \big( \max_{k \geq k_{s}} y_k -
	\min_{k \geq k_{s}} y_k\big)$ the signal amplitude in steady state, the relative
	RMS error is
\begin{equation} \label{eq:rom_rms}
	\epsilon_{\mathrm{rel}} = \frac{\big(\frac{1}{K_s} \sum_{k=k_s}^{K} \|y_{r,k} - y_k\|_2^2 \big)^{1/2}}{\frac{1}{2} \big( \max_{k \geq k_{s}} y_k - \min_{k \geq k_{s}} y_k\big)},
\end{equation}
where $K_s = K + 1 - k_s$ is the number of samples in the steady-state regime.

Three examples are considered. The first, in \Cref{ssec:cart}, serves as a proof
of concept on a low-dimensional nonlinear system for which an analytical solution
is known, demonstrating the algorithm's ability to recover exact solutions when the
components of $\pi$ belong to the span of the basis. The second and third examples,
in \Cref{ssec:lin} and \Cref{ssec:vdp}, consider a high-dimensional system of
dimension $n = 1000$ driven by a linear and a nonlinear signal generator
respectively, demonstrating scalability and the ability to construct accurate
reduced-order models in the absence of an analytical solution.

\subsection{Test 1: Cart pendulum}
\label{ssec:cart}

First, we consider the position control of a nonlinear cart pendulum, based on a
problem studied in \cite[Section 3.2]{romero_ii}. The signal generator is chosen
by selecting target dynamics for the cart pendulum with $d = 2$ and $n = 4$. The
dynamics of the interconnected system are governed by
\begin{equation*}
	\begin{aligned}
		s(\omega) = \begin{bmatrix}
			\omega_2 \\
			\frac{a_1 \sin(\omega_1)}{1 + k a_2 \cos(\omega_1)}
		\end{bmatrix}, && \ell(\omega) = \frac{k a_1 \sin(\omega_1)}{1 + k a_2 \cos(\omega_1)}, && k < - \frac{1}{a_2}
	\end{aligned}
\end{equation*}
\begin{equation*}
	\begin{aligned}
		f(x,u) = \begin{bmatrix}
			x_3 \\
			x_4 \\
			a_1 \sin(x_1) - a_2 \cos(x_1) u \\ 
			u
		\end{bmatrix} , && h(x) = x_1 &&  a_1 > 0, a_2 > 0 .
	\end{aligned}
\end{equation*}
This problem admits the analytical solution
\begin{equation*}
	\pi(\omega) = \left(\omega_1, k \omega_1, \omega_2, k \omega_2 \right)^\top .
\end{equation*}
The function $F$ and its Jacobian $JF$ take the form
\begin{equation*}
	\begin{aligned}
		F(\mathbf{c}) = \begin{bmatrix}
			\mathbf{A} \mathbf{c}_1 - \mathbf{M} \mathbf{c}_3 \\ \mathbf{A} \mathbf{c}_2 - \mathbf{M} \mathbf{c}_4 \\ \mathbf{A} \mathbf{c}_3 - H(\mathbf{c}_1) \\ \mathbf{A} \mathbf{c}_4 - \boldsymbol{\gamma}
		\end{bmatrix}, \quad 
		JF(\mathbf{c}) = \begin{bmatrix}
			\mathbf{A} & 0 & -\mathbf{M} & 0 \\
			0 & \mathbf{A} & 0 & -\mathbf{M} \\
			-\frac{\partial H}{\partial \mathbf{c}_1} & 0 & \mathbf{A} & 0 \\ 
			0 & 0 & 0 & \mathbf{A}
		\end{bmatrix}   
	\end{aligned}
\end{equation*}
using the previously defined matrix $\mathbf{A}$ in \eqref{eq:A}, as well as the
matrix and vector
\begin{equation} \label{eq:M_P_gamma}
	\begin{aligned}
		\mathbf{M} = \left \langle \Phi^N, \Phi^N \right \rangle_{\Omega}, && \boldsymbol{\gamma} = \left \langle \ell(\omega), \Phi^N \right \rangle_{\Omega} ,
	\end{aligned}
\end{equation}
and the function $H : \mathbb{R}^N \to \mathbb{R}^N$ given by
\begin{equation*}
	\begin{aligned}
		H(\mathbf{c}_1) & = \left \langle a_1 \sin \left(\Phi^N \mathbf{c}_1\right) - a_2 \cos\left(\Phi^N \mathbf{c}_1\right) \ell(\omega), \Phi^N \right \rangle _{\Omega} \\ 
		\frac{\partial H}{\partial \mathbf{c}_1} & = \left \langle \Phi^N \left(a_1 \cos \left(\Phi^N \mathbf{c}_1\right) + a_2 \sin\left(\Phi^N \mathbf{c}_1\right) \ell(\omega)\right), \Phi^N \right \rangle _{\Omega} .
	\end{aligned}
\end{equation*}
For this test case, the assembled matrix $\mathbf{A}$ and the resulting Jacobian $JF$ are both observed to be singular, and Newton iteration updates are therefore computed in the least-squares sense.
Using the constants $a_1 = 2$, $a_2 = 3$, and $k = -
\frac{2}{3}$, Table~\ref{table:cart} shows that the method achieves near-exact
recovery of the analytical solution over $\Omega = [-1,1]^2$. The solution
admits a monomial representation and the algorithm recovers it to within the
prescribed numerical tolerance, with observed errors on the order of $10^{-9}$.

\renewcommand{\arraystretch}{1.3}
\begin{table}[hbt!]
\caption{Residual error $\mathcal{E}_R(W)$ and
		computational time for the Galerkin approximation of the invariant mapping $\pi$
		for Test 1.}
	\begin{center}
		\begin{tabular}{ccccc}
			\hline
			&   $M = 2$ & $M = 4$ & $M = 6$ & $M = 8$ \\  \hline
			$\mathcal{E}_R(W)$   & $9.3 \times 10^{-10}$ & $8.6 \times 10^{-10}$ & $9.4 \times 10^{-10}$ & $2.4 \times 10^{-9}$ \\
			Time (s) & $0.83$ & $16.8$ & $135$ & $646$
		\end{tabular}
	\end{center}
	\label{table:cart}
\end{table}
\renewcommand{\arraystretch}{1}

\subsection{Test 2: Nonlinear resistor-inductor ladder with linear oscillator as signal generator} 
\label{ssec:lin}

We consider the nonlinear resistor-inductor (RL) ladder, based on a slightly
modified version of the circuit schematic in \cite[Figure 7]{RL_ladder} and
also considered in \cite{NL_MOR_kawano} in a model order reduction context. The
signal generator has dimension $d = 2$ and the system has a repeating structure
admitting any state dimension $n$. The signal generator is given by
\begin{equation}
	\begin{aligned}
		s(\omega) = \begin{bmatrix}
			0 & a \\ 
			-a & 0
		\end{bmatrix} \begin{bmatrix}
			\omega_1 \\ \omega_2
		\end{bmatrix}, & & \ell(\omega) = \omega_2, & & a \in \mathbb{R} \setminus\{ 0\}
	\end{aligned}
\end{equation}
and the system dynamics are
\begin{equation} \label{eq:RL_circuit}
	\begin{gathered}
		\begin{bmatrix}
			\dot{x}_1 \\ \dot{x}_2 \\ \vdots \\ \dot{x}_{n-1} \\ \dot{x}_n
		\end{bmatrix} = \underbrace{\begin{bmatrix}
				-2 \kappa & 1 \\ 
				1 & -2 \kappa & 1 \\
				& \ddots & \ddots & \ddots \\
				& & 1 & -2 \kappa & 1 \\
				& & & 1 & -2 \kappa
			\end{bmatrix} \begin{bmatrix}
				x_1 \\ x_2 \\ \vdots \\ x_{n-1} \\ x_n
			\end{bmatrix} - \begin{bmatrix}
				x_1^2/2 + x_1^3/3 \\ 
				x_2^2/2 + x_2^3/3 \\
				\vdots \\
				x_{n-1}^2/2 + x_{n-1}^3/3 \\
				x_{n}^2/2 + x_n^3/3
			\end{bmatrix} + \begin{bmatrix}
				1 \\ 0 \\ \vdots \\ 0 \\ 0
			\end{bmatrix} u}_\text{$=: f(x,u)$} \\
		h(x) = x_1 .
	\end{gathered}
\end{equation}
The parameter $\kappa = 1$ renders the system unstable beyond a certain dimension
$n$. Increasing this to $\kappa = 1.1$ was found to preserve stability for all values of $n$ considered in this study, up to $n = 1000$. We therefore use $\kappa = 1.1$ throughout this work. Since the nonlinearities in each
$f_i$ are polynomials of degree 2 and 3 depending only on $x_i$, the construction
of $F$ and $JF$ admits the following simplifications. We introduce the tensors
\begin{equation}
	\begin{aligned}
		\mathbf{N} \in \mathbb{R}^{N \times N \times N}, [\mathbf{N}]_{i,j,k} = \int_{\Omega} \phi_i \phi_j \phi_k d\omega \\ 
		\mathbf{O} \in \mathbb{R}^{N \times N \times N \times N}, [\mathbf{O}]_{i,j,k,l} = \int_{\Omega} \phi_i \phi_j \phi_k \phi_l d\omega ,
	\end{aligned}
\end{equation}
and the matrix-valued functions $\mathcal{N} : \mathbb{R}^N \to \mathbb{R}^{N \times N}$
and $\mathcal{O} : \mathbb{R}^N \to \mathbb{R}^{N \times N}$,
\begin{equation}
	\begin{aligned}
		\left[\mathcal{N}(\mathbf{v})\right]_{i,j} = \sum_{k=1}^N v_k N_{i,j,k}, &&& \left[\mathcal{O}(\mathbf{v})\right]_{i,j} = \sum_{k=1}^N \sum_{l=1}^N v_k v_l O_{i,j,k,l} .
	\end{aligned}
\end{equation}
Defining $P, Q : \mathbb{R}^N \to \mathbb{R}^{N \times N}$ as
\begin{equation}
	\begin{aligned}
		P(\mathbf{v}) & = \left(\mathbf{A} + 2 \kappa \mathbf{M} + \frac{1}{2} \mathcal{N}(\mathbf{v}) + \frac{1}{3} \mathcal{O}(\mathbf{v}) \right) \mathbf{v} \\
		Q(\mathbf{v}) & = \mathbf{A} + 2 \kappa \mathbf{M} + \mathcal{N}(\mathbf{v}) + \mathcal{O}(\mathbf{v}),   
	\end{aligned}
\end{equation}
where $\mathbf{A}$ and $\mathbf{M}$ are as in \eqref{eq:A} and \eqref{eq:M_P_gamma},
the function $F$ and its Jacobian take the form
\begin{equation} \label{eq:F_RL}
	F_i(\mathbf{c}) = P(\mathbf{c}_i) - (1 - \delta_{i,1}) \mathbf{M} \mathbf{c}_{i-1} - (1 - \delta_{i,n}) \mathbf{M} \mathbf{c}_{i+1} - \delta_{i,1} \boldsymbol{\gamma},
\end{equation}
where $\boldsymbol{\gamma}$ is the same as in \eqref{eq:M_P_gamma}, and $\delta_{i,j}$
	denotes the Kronecker delta, and
\begin{equation} \label{eq:JF_RL}
	JF(\mathbf{c}) = \begin{bmatrix}
		Q(\mathbf{c}_1) & -\mathbf{M} &&& \\[6pt]
		-\mathbf{M} & Q(\mathbf{c}_2) & -\mathbf{M} && \\[6pt]
		& \ddots & \ddots & \ddots & \\[6pt]
		&&  -\mathbf{M} & Q(\mathbf{c}_{n-1}) & -\mathbf{M} \\[6pt]
		&&& - \mathbf{M} & Q(\mathbf{c}_n) 
	\end{bmatrix}.
\end{equation}
Computing the tensors $\mathbf{N}$ and $\mathbf{O}$ exactly rather than by numerical
quadrature allows for a significant reduction in computational cost by exploiting
the structure of the system dynamics.

The method is validated on a full-order model with $n = 1000$ and
$a = 2$, over the domain $\Omega = [-0.6, 0.6]^2$, with residual
errors measured on $W = \Omega$. The choice of domain is discussed further
in \Cref{ssec:domain}. The residual errors $\mathcal{E}_R(W)$ reported in Table~\ref{table:rl_1000_tests} show that increasing $M$
systematically improves the approximation, as expected from the nested structure of
the polynomial spaces. The small residual errors at higher degrees suggest that
$[-0.6,0.6]^2$ lies within the neighbourhood $W^0$ over which $\pi$ is
well-defined, and that $\pi^N$ converges to the weak solution of the invariance
equation as the number of basis functions increases.

To construct the reduced-order model, we recall from \eqref{eq:constructed_rom}
the family of models
\begin{equation*} 
	\begin{aligned}
		\dot{r} & = s(r) - \bar{g}(r, \ell(r)) + \bar{g}(r, u) \\ 
		y_r & = h(\pi(r)) ,
	\end{aligned}
\end{equation*}
where $\bar{g}(\cdot)$ is a free mapping to be chosen to ensure local
exponential stability. For this system, the choice
\begin{equation}
	\bar{g}(r, u) = \begin{bmatrix}
		0 \\ c
	\end{bmatrix} u
\end{equation}
achieves local exponential stability for any $c > 0$. To verify this, we examine the linearization about the origin of $\frac{\partial}{\partial r} \left(s(r) - \begin{bmatrix} 0 \\ c \end{bmatrix} \ell(r)\right)$, given by the matrix $\begin{bmatrix} 0 & a \\ -a & -c \end{bmatrix}$, whose eigenvalues have strictly negative real part for all $c > 0$. The reduced-order model is then
\begin{equation}
	\begin{aligned}
		\dot{r} & = s(r) - \begin{bmatrix}
			0 \\ c
		\end{bmatrix} \ell(r) + \begin{bmatrix}
			0 \\ c
		\end{bmatrix} u \\
		y_r & = h(\pi^N(r ; \mathbf{c})),
	\end{aligned}
\end{equation}
for the approximate solution $\pi^N$ to the invariant mapping. Figure
\ref{fig:linear_rom} shows the outputs $y(t)$ and $y_r(t)$ for $n = 1000$, $M =
6$, $\Omega = [-0.6,0.6]^2$, $c = 10$, and initial conditions $\omega_0 = (0.2,
0.2)^\top$, $r_0 = (0,1)^\top$, $x_0 = 0$, over $t \in [0, 30]$. The
	relative RMS error $\epsilon_{\mathrm{rel}}$, defined in \eqref{eq:rom_rms} and reported in
	Table~\ref{table:rl_1000_tests}, is on the order of $10^{-3}$ across all values of
	$M$, confirming that the reduced-order model accurately reproduces the full-order
steady-state output despite the large dimensionality reduction.

\subsection{Test 3: Nonlinear resistor-inductor ladder with Van der Pol oscillator as signal generator}
\label{ssec:vdp}

We now consider the same RL ladder system \eqref{eq:RL_circuit} driven by the Van
der Pol oscillator as signal generator,
\begin{equation} \label{eq:vdp_s}
	\begin{aligned}
		s(\omega) = \begin{bmatrix}
			\omega_2 \\
			- \omega_1 + \mu (1 - \omega_1^2) \omega_2
		\end{bmatrix} & & \ell(\omega) = \omega_2 & & \mu > 0 .
	\end{aligned}
\end{equation}
The tensors $\mathbf{N}$ and $\mathbf{O}$ from Test 2 are reused; only $\mathbf{A}$ is
reconstructed using $s(\omega)$ in \eqref{eq:vdp_s}, and the expressions
\eqref{eq:F_RL}--\eqref{eq:JF_RL} for $F_i$ and $JF$ remain unchanged.

The invariance equation is solved with $\mu = 1$, $\kappa = 1.1$, over $\Omega =
W = [-0.6, 0.6]^2$ for $n = 1000$. The residual errors in
Table~\ref{table:rl_1000_tests} confirm that higher polynomial degree yields better
approximations, and that the method remains tractable under nonlinear excitation.

For the reduced-order model, we use the mapping
\begin{equation}
	\bar{g}(r, u) = \begin{bmatrix}
		0  \\ \mu (1 - r_1^2) + c
	\end{bmatrix} u,
\end{equation}
yielding the reduced dynamics
\begin{equation}
	\begin{aligned}
		\dot{r} & = s(r) - \begin{bmatrix}
			0  \\ \mu (1 - r_1^2) + c
		\end{bmatrix} \ell(r) + \begin{bmatrix}
			0  \\ \mu (1 - r_1^2) + c
		\end{bmatrix} u \\
		y_r & = h(\pi^N(r ; \mathbf{c})).
	\end{aligned}
\end{equation}
Local exponential stability follows from the linearization of
	$$\frac{\partial}{\partial r} \left(s(r) - \begin{bmatrix} 0 \\ \mu(1 - r_1^2) + c
	\end{bmatrix} \ell(r)\right)$$ about the origin, which gives the matrix
	$\begin{bmatrix} 0 & 1 \\ -1 & -c \end{bmatrix}$, whose eigenvalues have negative
	real part for all $c > 0$.

Figure~\ref{fig:vdp_rom} shows the model outputs for $n = 1000$, $M = 6$, $c =
10$, and initial conditions $\omega_0 = (0.2, 0.2)^\top$, $r_0 = (0,1)^\top$,
$x_0 = 0$, over $t \in [0, 30]$. The relative RMS errors $\epsilon_{\mathrm{rel}}$
in Table~\ref{table:rl_1000_tests} decrease monotonically with $M$, confirming
that the polynomial basis has sufficient representational capacity to capture the
nonlinear steady-state behaviour induced by the Van der Pol excitation.

\begin{table}[hbt!]
\caption{Galerkin approximation results for the RL ladder system ($n=1000$):
		residual error $\mathcal{E}_R(W)$ of the invariant mapping $\pi$
		and relative steady-state RMS error $\epsilon_{\mathrm{rel}}$ of the reduced-order
		model, for Test 2 (linear oscillator) and Test 3 (Van der Pol signal generator).}
	\centering
	\renewcommand{\arraystretch}{1.3}
	\begin{tabular}{llcccc}
		\hline
		& & $M=2$ & $M=4$ & $M=6$ & $M=8$ \\
		\hline
		\multirow{2}{*}{$\mathcal{E}_R(W)$}
		& Test 2 & $3.4\times10^{-4}$ & $1.9\times10^{-6}$ & $2.0\times10^{-8}$ & $6.5\times10^{-10}$ \\
		& Test 3 & $7.4\times10^{-3}$ & $2.0\times10^{-4}$ & $1.2\times10^{-5}$ & $8.1\times10^{-7}$ \\
		\hline
		\multirow{2}{*}{$\epsilon_{\mathrm{rel}}$}
		& Test 2 & $4.6\times10^{-3}$ & $4.5\times10^{-3}$ & $4.5\times10^{-3}$ & $4.5\times10^{-3}$ \\
		& Test 3 & $1.7\times10^{-1}$ & $7.0\times10^{-2}$ & $4.9\times10^{-2}$ & $4.1\times10^{-2}$ \\
		\hline
	\end{tabular}
	\label{table:rl_1000_tests}
\end{table}
\renewcommand{\arraystretch}{1}

\begin{figure}[htbp!]
	\centering
	\begin{subfigure}{0.8\textwidth}
		\centering
        \resizebox{\linewidth}{!}{\input{signal_plot_linear_hd.pgf}}
		\caption{Test 2: Linear oscillator signal generator.}
		\label{fig:linear_rom}
	\end{subfigure}
	\hfill
	\begin{subfigure}{0.8\textwidth}
		\centering
        \resizebox{\linewidth}{!}{\input{signal_plot_vdp_hd.pgf}}
		\caption{Test 3: Van der Pol signal generator.}
		\label{fig:vdp_rom}
	\end{subfigure}
	\caption{Time evolution of the excitation signal $\omega_1(t)$, the full-order
		output $y(t)$, the reduced-order output $y_r(t)$, and the error $|y(t)-y_r(t)|$
		for the RL ladder system ($n=1000$) over $t\in[0,30]$.}
	\label{fig:rom_comparison}
\end{figure}

\section{Discussion of Results} \label{sec:discussion}

The numerical results of the previous section demonstrate that the proposed method
achieves moment matching and accurately recovers the steady-state output of the
full-order model. The residual analysis confirms convergence of the polynomial
approximation to the weak solution of the invariance equation with increasing
degree. This section provides additional numerical evidence and examines the
behaviour of the method under parameter variations, with all experiments based on
the RL ladder driven by the Van der Pol oscillator from \Cref{ssec:vdp}, with the
same constants and initial conditions.

\subsection{Effect of domain and nonlinearities} \label{ssec:domain}

The mapping $\pi$ is only theoretically guaranteed to exist in a neighbourhood
$W^0$ of the origin, and the polynomial approximation $\pi^N$ is constructed over
a chosen domain $\Omega$. The size of $\Omega$ therefore directly affects both the
approximation quality of $\pi^N$ and whether the mapping satisfying
\eqref{eq:pde_sec3} is well-defined over that domain. Larger domains amplify
higher-order monomials near the boundary and can magnify small coefficient errors.

Figure~\ref{fig:domain_effect} illustrates this for $\mu = 0.25$, $\mu = 1$, and
$\mu = 2$, corresponding to increasingly nonlinear regimes of the Van der Pol
signal generator, with residuals measured on $W = [-0.6, 0.6]^2$. The results
confirm that the choice of $\Omega$ has a significant effect on approximation
quality, and that increased signal generator nonlinearity heightens sensitivity to
the domain. This is attributable to the cubic term $\mu \omega_1^2 \omega_2$
in the nonlinearity, which grows rapidly near the boundary when $\Omega$ or $\mu$
are large, thereby amplifying small coefficient errors in the approximate mapping $\pi^N$. Despite this, the method
maintains stable convergence and acceptable accuracy across the tested parameter
range for moderate domain sizes.

\begin{figure}[h]
	\centering
	\resizebox{0.8\textwidth}{!}{\input{mu_domain_hd.pgf}}
	\caption{PDE residual error as a function of domain size $\Omega$ for different
		values of the nonlinearity constant $\mu$ in the Van der Pol oscillator signal
		generator for Test 3. The tests were carried out with a maximum degree $M = 6$
		for the problem with $n = 1000$ dimensions and constants $\kappa = 1.1$, $c =
		10$, and tested on $W = [-0.6, 0.6]^2$.}
	\label{fig:domain_effect}
\end{figure}

\subsection{Choice of ansatz functions} \label{ssec:ansatz}

We note that the approximation quality of this method also depends on the choice of basis functions for the approximated invariant mapping $\pi^N$. In this work, monomial basis functions were selected primarily because many problems of interest exhibit polynomial nonlinearities, and a monomial basis naturally aligns with this structure, allowing for exact recovery of the corresponding polynomial solutions when they arise. In addition, their simple differentiation and integration properties allow the operators appearing in $F$ and $JF$ to be easily computed; for example, on symmetric domains such as $\Omega = [-L,L]^d$, many entries of $\mathbf{A}$ may vanish due to parity properties of the integrands for common signal generators $s$, promoting sparsity in these objects. We remark that, for a fixed total degree $M$, monomial and orthogonal polynomial bases span the same approximation space, so that in exact arithmetic the Galerkin approximation $\pi^N$ is independent of this choice. The basis affects only the conditioning of the resulting algebraic system, which is benign for the moderate degrees and domain sizes considered here. Alternative polynomial bases, such as Legendre or Chebyshev polynomials, could therefore offer improved numerical conditioning for higher degrees or larger domains, while non-polynomial ansatz functions, such as radial basis functions, would alter the approximation space itself and may be preferable when $\pi$ exhibits features poorly captured by globally supported polynomials. A systematic comparison of basis functions is beyond the scope of the present study.

\subsection{Comparison to other ROM techniques}

To provide a benchmark, the proposed method is compared against Proper Orthogonal
Decomposition (POD)~\cite{peraire} and Balanced Truncation
(BT)~\cite{scherpen1993balancing}, two established linear model reduction techniques
that do not exploit the structured interconnection or nonlinear excitation present
in the considered systems. The implementations are described as follows.

For POD, snapshots of the full-order state $x(t)$ were collected by simulating a single trajectory with an initial condition varying from the one tested, over the time interval $t \in [0,50]$ with initial conditions $\omega_0 = (0.1, 0.3)^\top$ and $x_0 = 0$. The snapshot matrix was assembled from 500 uniformly spaced time samples and the POD basis was extracted via a truncated singular value decomposition retaining $v = 2$ modes. The reduced dynamics were obtained by Galerkin projection of the full nonlinear system onto the POD subspace.

For BT, the method was applied to the system linearized about the origin. The controllability and observability Gramians were computed analytically by solving the associated Lyapunov equations, and the system was truncated to retain $v = 2$ states. The input and output matrices $B$ and $C$ were taken from the linearization of $f$ and $h$, respectively.

For Test 3 with $M = 8$ and reduced dimension $v = 2$, the proposed method
achieves $\epsilon_{\mathrm{rel}} = 4.1 \times 10^{-2}$, compared to $5.7 \times 10^{-2}$
for POD and $1.2 \times 10^{-1}$ for BT. Figure~\ref{fig:pod_comparison} further illustrates the advantage of moment
matching through phase portrait comparisons: the state trajectory reconstructed via
$\pi^N$ more closely follows that of the full-order model than either POD or BT at
the same reduced dimension.

\begin{figure}[htbp]
	\centering
	\begin{subfigure}{0.7\textwidth}
		\centering
		\resizebox{\linewidth}{!}{\input{phase_pod_hd.pgf}}
	\end{subfigure}
	\hfill
	\begin{subfigure}{0.7\textwidth}
		\centering
		\resizebox{\linewidth}{!}{\input{phase_bt_hd.pgf}}
	\end{subfigure}
	\caption{Phase portraits of the full-order model $(x_1, x_2)$ and approximate
		mapping $(\pi_1(\omega), \pi_2(\omega))$ for Test 3, contrasted with the POD
		ROM (top) and the BT ROM (bottom).}
	\label{fig:pod_comparison}
\end{figure}

\subsection{Computational complexity}

Table~\ref{table:computation} reports the number of basis functions $N$, Newton
iterations, and total computation time for Test 3. The Jacobian $JF$ is stored
sparsely, so the time complexity scales linearly in $n$ and quadratically in $N$.
The rapid convergence of both the residual and RMS error with $M$ means that $M =
4$ or $M = 6$ already provides a good approximation at manageable cost.

\renewcommand{\arraystretch}{1.3}
\begin{table}[hbt!]
\caption{Number of basis functions $N$, Newton iterations, and total computational
		time for the Galerkin approximation of $\pi$ in Test 3.}
	\begin{center}
		\begin{tabular}{ccccc}
			\hline
			& $M = 2$ & $M = 4$ & $M = 6$ & $M = 8$ \\ 
			\hline
			N  & $5$ & $14$ & $27$ & $44$ \\
			Iterations & 3 & 3 & 3 & 3 \\
			Time (s) & $1.92$ & $80$ & $1065$ & $8423$
		\end{tabular}
	\end{center}
	\label{table:computation}
\end{table}
\renewcommand{\arraystretch}{1}

Regarding online cost, for $t = 50$, $n = 1000$, and $v = 2$, the FOM requires
$0.4243$ seconds to simulate while the ROM requires $0.0120$ seconds, a speedup of
$35.3\times$. This arises from integrating a $v$-dimensional system rather than
the full $n$-dimensional one, with $v \ll n$, with the evaluation of $\pi^N$ at
the final time adding a cost that scales with $N$. For the offline
computation, the tensors $\mathbf{N}$ and $\mathbf{O}$ have sizes $O(N^3)$ and
$O(N^4)$ respectively; explicit storage was feasible for the problems considered
here, though symmetry-based compression or on-demand computation could be employed
if storage becomes prohibitive. The Jacobian $JF$ and function evaluation $F$ have
sizes $O(N^2)$ and $O(N)$ respectively, while the coefficients of $\pi^N$ require
$O(nN)$ storage online.

\subsection{Concluding remarks}

This paper presents a numerical framework for approximating the solution of the
invariance PDE system arising in nonlinear moment matching. A Galerkin expansion
through a polynomial basis yields an approximate invariant mapping $\pi^N$ that
can be used to construct reduced-order models accurately reproducing the
steady-state output of systems with state dimension up to $n = 1000$. Computational
efficiency is achieved by exploiting the problem structure through precomputed
matrix and tensor representations.

Nonlinear moment matching is most applicable when the dominant behaviour of a
high-dimensional system is governed by a low-dimensional driving subsystem. In
this setting, the reduced-order model evolves in a low-dimensional space while the
nonlinear complexity of the full-order system is encoded in the mapping $\pi$.
Despite the breadth of the theoretical moment matching literature, computational
methods for solving the invariance equation in high dimensions have been lacking;
this paper addresses that gap directly.

This moment matching method improves upon both POD and BT, and its underlying
invariance PDE formulation extends naturally beyond model reduction to output
regulation and immersion and invariance design. Future work will address signal
generators of higher dimensionality, enabling reduced-order models for richer input
regimes while requiring new strategies to manage the associated computational
challenges.

\section*{Code and Data Availability}

The code used to generate the numerical results in this paper is publicly available and can be viewed on GitHub at: \url{https://github.com/carlosjdoebeli/Moment-Matching-ROM}. 

\section*{Acknowledgments} D. Kalise is supported by the Air Force Office of Scientific Research under award number FA8655-26-1-B011. The authors would like to thank two anonymous reviewers, whose valuable feedback led to the final version of this manuscript.

\bibliographystyle{siamplain}
\bibliography{references}

@article{bungartz_griebel,
  author  = {Bungartz, Hans-Joachim and Griebel, Michael},
  title   = {Sparse grids},
  journal = {Acta Numerica},
  volume  = {13},
  pages   = {147--269},
  year    = {2004},
  doi     = {10.1017/S0962492904000182}
}

@article{dick_kuo_sloan,
  author  = {Dick, Josef and Kuo, Frances Y. and Sloan, Ian H.},
  title   = {High-dimensional integration: The quasi-{M}onte {C}arlo way},
  journal = {Acta Numerica},
  volume  = {22},
  pages   = {133--288},
  year    = {2013},
  doi     = {10.1017/S0962492913000044}
}

@article{kalise,
author = {Kalise, Dante and Kunisch, Karl},
year = {2017},
month = {02},
pages = {},
title = {Polynomial Approximation of High-Dimensional {Hamilton--Jacobi--Bellman} Equations and Applications to Feedback Control of Semilinear Parabolic {PDEs}},
volume = {40},
journal = {SIAM Journal on Scientific Computing},
doi = {10.1137/17M1116635}
}

@phdthesis{beard_thesis,
    author = {Randal W. Beard},
    title = {Improving the Closed-Loop Performance of Nonlinear Systems},
    school = {Rensselaer Polytechnic Institute},
    year = {1995}
}

@article{beard_galerkin,
title = {Galerkin approximations of the generalized {Hamilton-Jacobi-Bellman} equation},
author = {Randal W. Beard and George N. Saridis and John T. Wen},
journal = {Automatica},
volume = {33},
number = {12},
pages = {2159-2177},
year = {1997},
issn = {0005-1098},
doi = {10.1016/S0005-1098(97)00128-3},
}

@ARTICLE{astolfi_ii,
  author={Astolfi, A. and Ortega, R.},
  journal={IEEE Transactions on Automatic Control}, 
  title={Immersion and invariance: a new tool for stabilization and adaptive control of nonlinear systems}, 
  year={2003},
  volume={48},
  number={4},
  pages={590-606},
  doi={10.1109/TAC.2003.809820}}

@ARTICLE{astolfi_mm,
  author={Astolfi, Alessandro},
  journal={IEEE Transactions on Automatic Control}, 
  title={Model Reduction by Moment Matching for Linear and Nonlinear Systems}, 
  year={2010},
  volume={55},
  number={10},
  pages={2321-2336},
  doi={10.1109/TAC.2010.2046044}}

@INPROCEEDINGS{astolfi_mm_nonlinear,
  author={Astolfi, Alessandro},
  booktitle={2008 47th IEEE Conference on Decision and Control}, 
  title={Model reduction by moment matching for nonlinear systems}, 
  year={2008},
  volume={},
  number={},
  pages={4873-4878},
  doi={10.1109/CDC.2008.4738791}}

@book{carr,
  title     = "Applications of Centre Manifold Theory",
  author    = "Carr, Jack",
  year      = {1981},
  publisher = "Springer-Verlag New York Inc.",
  doi = {10.1007/978-1-4612-5929-9}
}

@book{large_scale_dynam,
author = {Antoulas, Athanasios C.},
title = {Approximation of Large-Scale Dynamical Systems},
publisher = {Society for Industrial and Applied Mathematics},
year = {2005},
doi = {10.1137/1.9780898718713}
}

@article{rafiq2022,
    author = {Rafiq, Danish and Bazaz, Mohammad Abid},
    title = {Model Order Reduction via Moment-Matching: A State of the Art Review},
    journal = {Archives of Computational Methods in Engineering},
    year = {2022},
    doi = {10.1007/s11831-021-09618-2}
}

@misc{kramer2022nonlinear,
      title={Nonlinear Balanced Truncation: Part 1-Computing Energy Functions}, 
      author={Boris Kramer and Serkan Gugercin and Jeff Borggaard and Linus Balicki},
      year={2022},
      eprint={2209.07645},
      archivePrefix={arXiv},
      primaryClass={math.OC}
}

@misc{kramer2023nonlinear,
      title={Nonlinear Balanced Truncation: Part 2 -- Model Reduction on Manifolds}, 
      author={Boris Kramer and Serkan Gugercin and Jeff Borggaard},
      year={2023},
      eprint={2302.02036},
      archivePrefix={arXiv},
      primaryClass={math.OC}
}

@article{antoulas_review,
title = {An overview of approximation methods for large-scale dynamical systems},
journal = {Annual Reviews in Control},
volume = {29},
number = {2},
pages = {181-190},
year = {2005},
issn = {1367-5788},
doi = {10.1016/j.arcontrol.2005.08.002},
author = {A.C. Antoulas},
}

@article{scarciotti2024interconnection,
  title={Interconnection-based model order reduction -- a survey},
  author={Scarciotti, Giordano and Astolfi, Alessandro},
  journal={European Journal of Control},
  volume={75},
  pages={100929},
  year={2024},
  publisher={Elsevier}
}

@book{isidori1995nonlinear,
author = {Isidori, Alberto},
title = {Nonlinear Control Systems},
publisher = {Springer London},
year = {1995},
doi = {10.1007/978-1-84628-615-5},
edition   = {3}
}

@article{romero_ii,
author = {Ortega, Romeo and Yi, Bowen and Romero, Jose Guadalupe and Astolfi, Alessandro},
title = {Orbital stabilization of nonlinear systems via the immersion and invariance technique},
journal = {International Journal of Robust and Nonlinear Control},
volume = {30},
number = {5},
pages = {1850-1871},
doi = {10.1002/rnc.4861},
year = {2020}
}

@ARTICLE{RL_ladder,
  author={Rewienski, M. and White, J.},
  journal={IEEE Transactions on Computer-Aided Design of Integrated Circuits and Systems}, 
  title={A trajectory piecewise-linear approach to model order reduction and fast simulation of nonlinear circuits and micromachined devices}, 
  year={2003},
  volume={22},
  number={2},
  pages={155-170},
  doi={10.1109/TCAD.2002.806601}}

@article{NL_MOR_kawano,
title = {Empirical differential Gramians for nonlinear model reduction},
journal = {Automatica},
volume = {127},
pages = {109534},
year = {2021},
issn = {0005-1098},
doi = {10.1016/j.automatica.2021.109534},
author = {Yu Kawano and Jacquelien M.A. Scherpen},
}

@ARTICLE{gray_energy,
  author={Scherpen, J.M.A. and Gray, W.S.},
  journal={IEEE Transactions on Automatic Control}, 
  title={Minimality and local state decompositions of a nonlinear state space realization using energy functions}, 
  year={2000},
  volume={45},
  number={11},
  pages={2079-2086},
  doi={10.1109/9.887630}}

@INPROCEEDINGS{gray_hilbert,
  author={Gray, W.S. and Scherpen, J.M.A.},
  booktitle={Proceedings of the 2001 American Control Conference.}, 
  title={Nonlinear {Hilbert} adjoints: properties and applications to {Hankel} singular value analysis}, 
  year={2001},
  volume={5},
  number={},
  pages={3582-3587},
  doi={10.1109/ACC.2001.946190}}

@article{gray_manifolds,
author = {Gray, W. and Verriest, Erik},
year = {2006},
month = {04},
pages = {653-659},
title = {Balanced Realizations Near Stable Invariant Manifolds},
volume = {42},
journal = {Automatica},
doi = {10.1016/j.automatica.2005.12.007}
}

@article{peraire,
author = {Willcox, K. and Peraire, J.},
year = {2002},
month = {11},
pages = {2323-2330},
title = {Balanced Model Reduction via the Proper Orthogonal Decomposition},
volume = {40},
journal = {AIAA Journal},
doi = {10.2514/2.1570}
}

@article{Kunisch1999ControlOT,
  title={Control of the {B}urgers equation by a reduced-order approach using proper orthogonal decomposition},
  author={Karl Kunisch and Stefan Volkwein},
  journal={Journal of Optimization Theory and Applications},
  year={1999},
  volume={102},
  pages={345-371},
  url={https://api.semanticscholar.org/CorpusID:118949739}
}

@article{Kunisch2008,
author = {Karl Kunisch and Stefan Volkwein},
journal = {{ESAIM}: Mathematical Modelling and Numerical Analysis},
language = {eng},
month = {1},
number = {1},
pages = {1-23},
publisher = {EDP Sciences},
title = {Proper orthogonal decomposition for optimality systems},
url = {http://eudml.org/doc/250374},
volume = {42},
year = {2008},
}

@book{khalil,
  title={Nonlinear Systems},
  author={Khalil, H.K.},
  isbn={9789332542037},
  series={Always Learning},
  url={https://books.google.co.uk/books?id=Gt2HAQAACAAJ},
  year={2015},
  publisher={Pearson}
}

@book{bellman,
 ISBN = {9780691079011},
 URL = {http://www.jstor.org/stable/j.ctt183ph6v},
 author = {Richard Bellman},
 publisher = {Princeton University Press},
 title = {Adaptive Control Processes: A Guided Tour},
 urldate = {2023-10-27},
 year = {1961}
}

@book{MOR_app,
  title={Model Order Reduction: Volume 3 Applications},
  author={Benner, Peter and Schilders, Wil and Grivet-Talocia, Stefano and Quarteroni, Alfio and Rozza, Gianluigi and Miguel Silveira, Lu{\'\i}s},
  year={2020},
  publisher={De Gruyter},
isbn = {9783110499001},
doi = {10.1515/9783110499001}
}

@article{ionescu2016nonlinear,
    title={Nonlinear moment matching-based model order reduction},
    author={Ionescu, T. C. and Astolfi, A.},
    journal={IEEE Transactions on Automatic Control},
    volume={61},
    number={10},
    pages={2837--2847},
    year={2015},
    publisher={IEEE}
}

@article{scarciotti2017data,
title = {Data-driven model reduction by moment matching for linear and nonlinear systems},
journal = {Automatica},
volume = {79},
pages = {340-351},
year = {2017},
issn = {0005-1098},
doi = {https://doi.org/10.1016/j.automatica.2017.01.014},
author = {Giordano Scarciotti and Alessandro Astolfi}
}

@article{scherpen1993balancing,
	title={Balancing for nonlinear systems},
	author={Scherpen, J. M. A.},
	journal={Systems \& Control Letters},
	volume={21},
	number={2},
	pages={143--153},
	year={1993},
	publisher={Elsevier}
}

@inproceedings{simard2023parameterization,
	title={Parameterization of All Moment Matching Interpolants},
	author={Simard, J. D. and Moreschini, A. and Astolfi, A.},
	booktitle={European Control Conference (ECC)},
	pages={1--6},
	year={2023},
  doi={10.23919/ECC57647.2023.10178406}
}

@article{simard2023parameterizationb,
  author={Simard, Joel D. and Moreschini, Alessio and Astolfi, Alessandro},
  journal={IEEE Transactions on Automatic Control}, 
  title={Parameterization of All Differential-Algebraic Moment Matching Interpolants}, 
  year={2025},
  volume={70},
  number={3},
  pages={1875-1882},
  doi={10.1109/TAC.2024.3469247},
}

@article{moreschini2024closed,
 author={Moreschini, A. and Astolfi, A.},
  journal={IEEE Transactions on
Automatic Control}, 
  title={Closed-Loop Interpolation by Moment Matching for Linear and Nonlinear Systems}, 
year={2025},
  volume={70},
  number={5},
  pages={2918-2933},
  doi={10.1109/TAC.2024.3484309}
}

@ARTICLE{moreschini2024a,
  author={Moreschini, Alessio and Scandella, Matteo and Astolfi, Alessandro and Parisini, Thomas},
  journal={IEEE Transactions on Automatic Control}, 
  title={Moment Matching by {Kernel}-Based Learning}, 
  year={2026},
  volume={71},
  number={4},
  pages={2123-2138},
  doi={10.1109/TAC.2025.3618165}
  }

@misc{kalise2019robust,
      title={Robust feedback control of nonlinear {PDEs} by numerical approximation of high-dimensional {Hamilton-Jacobi-Isaacs} equations}, 
      author={Dante Kalise and Sudeep Kundu and Karl Kunisch},
      year={2019},
      eprint={1905.06276},
      archivePrefix={arXiv},
      primaryClass={math.OC}
}

@article{byrnes2003limit,
  title={Limit sets, zero dynamics, and internal models in the problem of nonlinear output regulation},
  author={Byrnes, Christopher I and Isidori, Alberto},
  journal={IEEE Transactions on Automatic Control},
  volume={48},
  number={10},
  pages={1712--1723},
  year={2003},
  publisher={IEEE}
}

@article{isidori2008steady,
  title={Steady-state behaviors in nonlinear systems with an application to robust disturbance rejection},
  author={Isidori, Alberto and Byrnes, Christopher I},
  journal={Annual Reviews in Control},
  volume={32},
  number={1},
  pages={1--16},
  year={2008},
  publisher={Elsevier},
 doi={10.1016/j.arcontrol.2008.01.001}
}

@article{faedo_scarciotti_astolfi,
  author={Faedo, Nicolás and Scarciotti, Giordano and Astolfi, Alessandro and Ringwood, John V.},
  journal={IEEE Transactions on Automatic Control}, 
  title={On the Approximation of Moments for Nonlinear Systems}, 
  year={2021},
  volume={66},
  number={11},
  pages={5538-5545},
  doi={10.1109/TAC.2021.3054325}
}

@misc{benner,
      title={Interpolation-Based Model Order Reduction for Polynomial Parametric Systems}, 
      author={Peter Benner and Pawan Goyal},
      year={2019},
      eprint={1904.11891},
      archivePrefix={arXiv},
      primaryClass={math.NA} 
}

@article{gallivan,
title = {Sylvester equations and projection-based model reduction},
journal = {Journal of Computational and Applied Mathematics},
volume = {162},
number = {1},
pages = {213-229},
year = {2004},
note = {Proceedings of the International Conference on Linear Algebra and Arithmetic 2001},
issn = {0377-0427},
doi = {10.1016/j.cam.2003.08.026},
author = {K. Gallivan and A. Vandendorpe and P. {Van Dooren}}
}

@ARTICLE{hermann,
  author={Hermann, R. and Krener, A.},
  journal={IEEE Transactions on Automatic Control}, 
  title={Nonlinear controllability and observability}, 
  year={1977},
  volume={22},
  number={5},
  pages={728-740},
  doi={10.1109/TAC.1977.1101601}}

@article{sussmann,
title = {Controllability of nonlinear systems},
journal = {Journal of Differential Equations},
volume = {12},
number = {1},
pages = {95-116},
year = {1972},
issn = {0022-0396},
doi = {10.1016/0022-0396(72)90007-1},
url = {https://www.sciencedirect.com/science/article/pii/0022039672900071},
author = {H\'ector J Sussmann and Velimir Jurdjevic}
}

@article{AKK21,
  author  = {Behzad Azmi and Dante Kalise and Karl Kunisch},
  title   = {Optimal Feedback Law Recovery by Gradient-Augmented Sparse Polynomial Regression},
  journal = {Journal of Machine Learning Research},
  year    = {2021},
  volume  = {22},
  number  = {48},
  pages   = {1--32},
  url     = {http://jmlr.org/papers/v22/20-755.html}
}
\end{document}